\documentclass[12pt]{article}

\title{Extended real number arithmetics via Dedekind cuts}

\author{
Andreas H. Hamel\footnote{Free University of Bolzano-Bozen, Faculty of Economics and Management, \href{mailto:andreas.hamel@unibz.it}{andreas.hamel@unibz.it}}
} 
\date{{\small \today}}

\usepackage{array} % for better arrays (eg matrices) in maths
\usepackage{verbatim} % adds environment for commenting out blocks of text & for better verbatim

\usepackage{amsmath, amssymb, enumerate, color}
\usepackage{stmaryrd}

\newtheorem{theorem}{Theorem}
\newtheorem{corollary}[theorem]{Corollary}
\newtheorem{remark}[theorem]{Remark}

\newtheorem{definition}[theorem]{Definition}
\newtheorem{proposition}[theorem]{Proposition}
\newtheorem{example}[theorem]{Example}

\numberwithin{equation}{section}  % Formeln mit f¸hrender Kapitelnummer
\numberwithin{figure}{section}    % Abbildungen mit f¸hrender Kapitelnummer
\numberwithin{table}{section}     % Tabellen mit f¸hrender Kapitelnummer
\numberwithin{theorem}{section}

\newcommand{\of}[1]{\ensuremath{\left( #1 \right)}}

\newcommand{\cb}[1]{\ensuremath{ \left\{ #1 \right\} }}

\newcommand{\bs}{\backslash}
\newcommand{\vn}{\varnothing}

\newcommand{\pend}{ \hfill $\square$ \medskip}

\newcommand{\eps}{\ensuremath{\varepsilon}}

\newcommand{\vp}{\ensuremath{\varphi}}

\renewcommand{\P}{\ensuremath{\mathcal{P}}}

\newcommand{\R}{\mathrm{I\negthinspace R}}

\newcommand{\dom}{{\rm dom \,}}
\newcommand{\epi}{{\rm epi \,}}

\newcommand{\cl}{{\rm cl \,}}

\newcommand{\co}{{\rm conv \,}}

\newcommand{\isum}{{+^{\negmedspace\centerdot\,}}}
\newcommand{\ssum}{{+_{\negmedspace\centerdot\,}}}
\newcommand{\idif}{{-^{\negmedspace\centerdot\,}}}
\newcommand{\sdif}{{-_{\negmedspace\centerdot\,}}}

\newcommand{\triup}{{\rm \vartriangle}}

\usepackage[%backref=page,
colorlinks=true, linkcolor=blue, citecolor=blue, urlcolor=blue, filecolor=blue
auskommentieren und unten pdfborder aktivieren]{hyperref}
\definecolor{color0}{gray}{.50}
\definecolor{color1}{rgb}{0,.2,.8}
\definecolor{color2}{rgb}{1,.2,0}
\definecolor{color3}{rgb}{.2,.7,.6}

\begin{document}

\maketitle

\begin{abstract}
It is shown how Dedekind cuts can be used to introduce the extended real numbers along with sound arithmetic laws via one simple rule for the addition of sets. The crucial idea is that the use of the lower and the upper part of the cuts, respectively, leads to two different additions which are known in the literature as inf-addition and sup-addition. Moreover, the two resulting structures are conlinear spaces which at the same time are complete lattices with respect to the natural order. This admits the definition of pseudo-differences on the extended reals which also provide formulas for expressions like $(+\infty) - (+\infty)$, $(-\infty) - (-\infty)$. There are two major motivations: one is that proper and improper extended real-valued functions can be treated in a unified manner, the other that set-valued functions can often be represented by families of scalar functions which may include improper ones. 
\end{abstract}



%%%New section
\section{Introduction}

Despite the fact that the extension of the real numbers by a bottom and a top element, usually denoted by $-\infty$ and $+\infty$, is extremely useful and popular in Variational Analysis and Optimization Theory, most authors have had a difficult time to introduce and justify proper arithmetic rules for the extended real number system. A few examples follow.

In \cite[p. 29]{AliprantisBurkinshaw98Book} one can read: `The expressions $\infty - \infty$ and $-\infty + \infty$ are left (as usual) undefined. In this book we shall agree that
\[
\negthickspace\negthickspace\negthickspace\negthickspace\negthickspace\negthickspace\negthickspace   5. \qquad 0 \cdot \infty = 0.'
\]
In \cite[p. 39]{Zalinescu02Book}, the definition of a convex function via Jensen's inequality is given `with the conventions: $(+\infty) +(-\infty) = +\infty$, $0 \cdot (+\infty) = +\infty$, $0 \cdot (-\infty) = 0$.'

Two well-known texts, Rockafellar's \cite[p. 6]{Rockafellar74Book} and Aubin's \cite[p. 10]{Aubin79BookRev} considered the arithmetic operation $(+\infty) + (-\infty)$ as 'risky' (Rockafellar) and even 'undefined and {\em forbidden}' (Aubin, emphasis in the original). 

Sometimes, the difficulty is just avoided. Phelps \cite[p. 40]{Phelps89Book} comforts his readers: 'we won't have to worry about $+\infty - \infty$. In \cite[p. 24]{Rockafellar70Book}, Rockafellar also writes 'The combinations $\infty - \infty$ and $- \infty + \infty$ are undefined and are avoided.' However, in the later book \cite[p. 6]{Rockafellar74Book} one can already find the comment that the rule $(+\infty) + (-\infty) = +\infty$, is suited for convex functions while $(+\infty) + (-\infty) = -\infty$, 'is needed for concave functions.' 

In the more recent book \cite[p. 15]{RockafellarWets10Book3rd} by Rockafellar and Wets, both versions for $(+\infty) + (-\infty)$ are discussed, but it is stated that 'there's no single, symmetric way of handling $\infty - \infty$.' The first variant, called inf-addition, is introduced as a 'convention' and mostly used in the book. A warning is given to the readers: 'The opposite convention is {\em sup-addition} (i.e., $(+\infty) + (-\infty) = -\infty$, emphasis in the original -- A.H.); we won't invoke it without explicit warning.'

The book \cite{AliprantisBorder06Book3rd} by Aliprantis and Border is very popular and widely used as a reference text. Therein, one can read already on page 2 'The combination $+\infty - \infty$ of symbols has no meaning. The symbols $+\infty$ and $-\infty$ are not really meant to be used for arithmetic, they are only used to avoid awkward expressions involving infima and suprema.'

Moreover, it has been realized that the extension of Dedekind cuts by the two "non-ordinary" cuts $(\mathbb Q, \vn)$ and $(\vn, \mathbb Q)$ might be required to resolve the issue but, to the best of our knowledge, no systematic attempt has been made to derive appropriate arithmetic rules for the extended real numbers based on extended Dedekind cuts. This gap is closed in this paper by giving a single set addition rule upon which two arithmetic structures for the extended reals can be based. The crucial observation is that one has to focus either on the "lower" or the "upper" half of the cut--which does not make a difference for ordinary Dedekind cuts, but it does for non-ordinary ones.

In a number of papers, the twofold extended arithmetic is already used. Moreau \cite{Moreau63SMUM, Moreau63RCAS, Moreau66SJL} was the first who applied an extended real number arithmetic to problems in Variational Analysis, especially infimal convolution. The observation that inf-addition is suited for convexity while sup-addition for concavity made in \cite{Rockafellar74Book, RockafellarWets10Book3rd} seems to be widely accepted in Variational Analysis. Further examples include the papers \cite{GetanMartinezLegazSinger03JMS} and \cite{CohenGaubertQuadratSinger05InColl} (with different purposes) while in \cite{HamelSchrage12JCA} much of the arithmetic was developed without any reference to Dedekind cuts or another construction method for the (extended) real numbers. A framework for algebraic structures which are at the same time complete lattices was developed in \cite{Hamel05Habil} motivated by applications in set optimization.

%%%New section
\section{Power semigroups}

Let $(W, +)$ be a commutative semigroup. The operation $+ \colon W \times W \to W$ is dubbed "addition." Commutativity is assumed since the main application, the extended Dedekind cuts, enjoy this property. The addition can be extended to the power set of $W$, denoted $\P(W)$, by
\begin{equation}
\label{EqMinkowskiAddition}
A + B = \cb{a + b \mid a \in A, \, b \in B}
\end{equation}
for $A, B \in \P(W)\bs\{\vn\}$ and
\begin{equation}
\label{EqEmptyAddition}
A + \vn = \vn + A = \vn
\end{equation}
for all $A \in \P(W)$. Especially, $\vn + \vn = \vn$. Equation \eqref{EqEmptyAddition} means that $\vn$ is the absorbing element in $(\P(W), +)$ which is unique and an idempotent (see, for example, \cite{HagerHamelHeyde23ArX} and the references therein).

If $W$ includes a neutral element with respect to addition, i.e., $0 \in W$ with $0 + 0 = 0$, then $\{0\} + \{0\} = \{0\}$, i.e., $\{0\}$ is the neutral element in $(\P(W), +)$. In this case, both $(W, +)$ and $(\P(W), +)$ are (commutative) monoids. This fact is well-known and stated for the sake of reference. The straightforward proof is omitted.

\begin{proposition}
\label{PropPowerSemigroup}
If $(W, +)$ is a (commutative) semigroup, then so is $(\P(W), +)$. If $(W, +)$ is a (commutative) monoid, then so is $(\P(W), +)$.
\end{proposition}

Equation \eqref{EqEmptyAddition} is the crucial agreement for the extension of Dedekind cuts from those involving only non-empty sets to ones with an empty lower or upper set. This will be done in the next section.

%%%New section
\section{Ordinary and non-ordinary Dedekind cuts}

Let $\mathbb Q$ denote the set of rational numbers. The usual order relations on $\mathbb Q$ are denoted $<$ (strictly less; asymmetric, transitive and total) and $\leq$ (less than or equal to; reflexive, transitive, antisymmetric and a linear (complete) order).

\begin{definition}
\label{DefDedekindCut}
A Dedekind cut is a pair $(A, B) \in \P(\mathbb Q) \times \P(\mathbb Q)$ satisfying

(D1) $q \in \mathbb Q$ implies $q \in A$ or $q \in B$ but not both,

(D2) $p \in A$, $q \in B$ implies $p < q$

and either

(DL3) if $a \in A$, then there is $a' \in A$ such that $a < a'$ ($A$ does not include a greatest element),

or 

(DU3) if $b \in B$, then there is $b' \in B$ such that $b' < b$ ($B$ does not include a least element).

A Dedekind cut $(A, B)$ is called ordinary if $A, B \not\in \{\vn, \mathbb Q\}$.  
\end{definition}

The collection of all Dedekind cuts satisfying (DL3) is denoted by $\mathcal D_L$, the ones with (DU3) by $\mathcal D_U$. There are exactly two non-ordinary Dedeking cuts, namely $(\vn, \mathbb Q)$ and $(\mathbb Q, \vn)$.

If $(A, B)$ is a Dedekind cut with either $A = \{p \in \mathbb Q \mid p \leq a\}$ for some $a \in \mathbb Q$ or $B = \{q \in \mathbb Q \mid b \leq q\}$ for some $b \in \mathbb Q$, then it is called generated by $a$ and $b$, respectively. 

A Dedekind cut can satisfy both of (DL3), (DU3), and these are exactly the ones which are not generated by rational numbers. In  particular, non-ordinary Dededind cuts are of this type. Let $(A, B)$ be an ordinary Dedekind cut which satisfies (DL3), but not (DU3). Then $B$ has a least element, say $\bar q$, i.e., $B = \{q \in \mathbb Q \mid \bar q \leq q\}$ and $(A \cup\{\bar q\}, B\bs\{\bar q\})$ is a Dedekind cut which satisfies (DU3), but not (DL3) which is generated by the same rational number $\bar q$.

To re-install uniqueness, it is customary to move from cuts $(A, B)$ to ``lower sets'' $A$ with (DL3) (done, e.g.,  in \cite{Rudin76Book}) or ``upper sets'' $B$ with (DU3) (done, e.g., in \cite{EbbinghausEtAl88Book2nd}). Either one of these two moves is enough to construct the real numbers as ordinary Dedekind cuts. It turns out that this produces two different results if non-ordinary cuts are taken into consideration.

\begin{proposition}
\label{PropLowerSets}
The pair $(A, B) \in \P(\mathbb Q) \times \P(\mathbb Q)$ is a Dedekind cut satisfying (DL3) if, and only if, 

(A1) $B = \mathbb Q\bs A$ 

(A2) $a, a' \in \mathbb Q$, $a' < a$ and $a \in A$ imply $a' \in A$ ($A$ is closed downward)

(A3) $A$ has no greatest element.

Such a pair is an ordinary Dedekind cut if, and only if, $A \not\in \{\vn, \mathbb Q\}$.
\end{proposition}

{\sc Proof.} If $A = \mathbb Q$, $B =\vn = \mathbb Q\bs A$, then $A$ satisfies (A1), (A2) as well as  (D2), (D3) and (DL4). The same is true for $A = \vn$, $B =\mathbb Q = \mathbb Q\bs A$. If $(A, B)$ is an ordinary Dedekind cut, then (A1) is immediate from (D1) and (A2) follows from (D1) and (D2) while (A3) is identical with (DL3). 

Vice versa, let $(A, B) \in \P(\mathbb Q) \times \P(\mathbb Q)$ satisfy (A1)-(A3). Then (D1) follows from (A1) and (D2) from (A1) and (A2) while (DL3) is identical with  (A3). 

The last claim is equally obvious. \pend

The collection of sets $A \subseteq \mathbb Q$ satisfying (A2) and (A3) in Proposition \ref{PropLowerSets} is denoted $\mathcal L$. Condition (A2) can also be stated as: $A$ is an element of 
\[
\mathcal P(\mathbb Q, -\mathbb Q_+) = \cb{D \subseteq \mathbb Q \mid D - \mathbb Q_+ \subseteq D}
\]
where the notation is consistent with the one used in \cite{HamelEtAl15Incoll} and $D - \mathbb Q_+ = D + (-\mathbb Q_+)$ with $-\mathbb Q_+ = \{p \in \mathbb Q \mid p \leq 0\}$. Elements of $\mathcal P(\mathbb Q, -\mathbb Q_+)$ may satisfy (A3) or not, i.e., $\mathcal L \subsetneq \mathcal P(\mathbb Q, -\mathbb Q_+)$. Therefore, the function $I \colon \mathcal P(\mathbb Q, -\mathbb Q_+) \to \mathcal L$ is introduced by
\[
I(A) = \left\{
	\begin{array}{cc}
	A & \text{if $A$ satisfies (A3)}, \\
	A\bs\{a\} & \text{if $a$ is the greatest element of $A$}.
	\end{array}
	\right.
\]
Especially, $I(\vn) = \vn$ and $I(\mathbb Q) = \mathbb Q$.

\begin{proposition}
\label{PropIntensiveClosure} The function $I \colon \mathcal P(\mathbb Q, -\mathbb Q_+) \to \mathcal L$ is a Kuratowski interior operator.
\end{proposition}

{\sc Proof.} Indeed, $I$ maps into $\mathcal L$. Moreover, it clearly satisfies $I(\mathbb Q) = \mathbb Q$, is intensive ($I(A) \subseteq A$ for all $A \in \mathcal P(\mathbb Q, -\mathbb Q_+)$), idempotent ($I(I(A)) = I(A)$ for all $A \in \mathcal P(\mathbb Q, -\mathbb Q_+)$) and preserves binary intersections, i.e., $A, A' \in \mathcal P(\mathbb Q, -\mathbb Q_+)$ implies $I(A \cap A') = I(A) \cap I(A')$: if $p \in I(A \cap A')$, then there is $q \in A \cap A'$ with $p < q$, hence $p \in I(A) \cap I(A')$, and, vice versa, if $p \in I(A) \cap I(A')$, then there are $a \in A$, $a' \in A'$ with $p < a$, $p < a'$ and hence $p < \min\{a, a'\} \in A \cap A'$ which is true since $A, A' \in \mathcal P(\mathbb Q, -\mathbb Q_+)$.

Thus, $I$ is an interior (or kernel) operator. \pend

Note that the preservation of binary intersections gives monotonicity ($A \subseteq A'$ implies $I(A) \subseteq I(A')$) which also follows from the very definition of $I$.

\begin{remark}
\label{RemUpperSets} The corresponding construction of ``upper sets" leads to the collection $\mathcal U \subseteq \mathcal P(\mathbb Q, \mathbb Q_+)$ whose elements are (B2) closed upward and (B3) have no least element. A counterpart for Proposition \ref{PropLowerSets} is straightforward. The corresponding interior operator  is $J \colon \mathcal P(\mathbb Q, \mathbb Q_+) \to \mathcal U$ defined by
\[
J(B) = \left\{
	\begin{array}{cc}
	B & \text{if $B$ satisfies (B3)}, \\
	B\bs\{b\} & \text{if $b$ is the least element of $B$}.
	\end{array}
	\right.
\]
\end{remark}

It is sufficient (and even necessary) to consider either $\mathcal L$ or $\mathcal U$ for the construction of real numbers as ordinary cuts. Indeed, the embedding of rational numbers into the set of cuts yields a unique one only if one considers cuts with (DL3) or with (DU3): $q \in \mathbb Q$ generates the cut $(\{p \in \mathbb Q \mid p < q\}, \{p \in \mathbb Q \mid q \leq p\})$ which satisfies (DL3) and $(\{p \in \mathbb Q \mid p \leq q\}, \{p \in \mathbb Q \mid q < p\})$ which satisfies (DU3). 

There is a one-to-one correspondence between cuts satisfying (DL3), elements of $\mathcal U$, cuts satisfying (DU3) and elements of $\mathcal L$.

\begin{theorem}
\label{ThmOneToOneStructures}
Any two of the following sets are in one-to-one correspondence: $\mathcal D_L$, $\mathcal D_U$, $\mathcal L$ and $\mathcal U$.
\end{theorem}

{\sc Proof.} Proposition \ref{PropLowerSets} and its ``upper" counterpart (see Remark \ref{RemUpperSets}) provide one-to-one correspondences between $\mathcal D_L$, $\mathcal L$ and $\mathcal D_U$, $\mathcal U$, respectively. The one-to-one correspondence between $\mathcal L$ and $\mathcal U$ is provided by the following formulas:
\begin{align*}
A \in \mathcal L & \Rightarrow J(\mathbb Q \bs A) \in \mathcal U \quad \text{and} \quad I\of{\mathbb Q \bs J(\mathbb Q \bs A)} = A, \\
B \in \mathcal U & \Rightarrow I(\mathbb Q \bs B) \in \mathcal L \quad \text{and} \quad J\of{\mathbb Q \bs I(\mathbb Q \bs B)} = B.
\end{align*}
Indeed, on the one hand $A = \mathbb Q \bs (\mathbb Q \bs A) \subseteq \mathbb Q \bs J(\mathbb Q \bs A)$ since $J(\mathbb Q \bs A) \subseteq \mathbb Q \bs A$. This gives $A = I(A) \subseteq I(\mathbb Q \bs J(\mathbb Q \bs A))$ for $A \in \mathcal L$. One the other hand, $a \in I(\mathbb Q \bs J(\mathbb Q \bs A))$ means $a \in \mathbb Q \bs J(\mathbb Q \bs A)$ and $a$ is not the greatest element of $\mathbb Q \bs J(\mathbb Q \bs A)$. Two cases are possible. First, $\mathbb Q \bs A$ does not have a least element in which case $J(\mathbb Q \bs A) = \mathbb Q \bs A$ which in turn gives $a \in \mathbb Q \bs (\mathbb Q \bs A) = A$. Secondly, $\mathbb Q \bs A$ has a least element $b \in \mathbb Q$ in which case $J(\mathbb Q \bs A) = (\mathbb Q \bs A) \bs \{b\}$. Since $b$ is the least element of $\mathbb Q \bs A$, it is the greatest element of $\mathbb Q \bs J(\mathbb Q \bs A) = A \cup\{b\}$, hence $a = b$ is not possible. Therefore, $a \in A$ also in this case. \pend

Following Dedekind \cite{Dedekind1872}, $\mathcal D_L\bs\{(\vn, \mathbb Q), (\mathbb Q, \vn)\}$ or $\mathcal D_U\bs\{(\mathbb Q, \vn), (\vn, \mathbb Q)\}$ are identified with the set of real numbers, see \cite{EbbinghausEtAl88Book2nd, Rudin76Book}. Because of Theorem \ref{ThmOneToOneStructures}, $\mathcal L\bs\{\mathbb Q,\vn\}$ or $\mathcal U\bs\{\mathbb Q,\vn\}$ can also be used. The addition of $\mathbb Q,\vn$ provides two more elements to be dealt with which are customarily denoted by $-\infty$, $+\infty$. This results in structures which have better order-type properties, i.e., they are complete lattices instead of what is nowadays called Dedekind (or conditionally) complete lattices and restricted algebraic features, i.e., they are only commutative monoid, not groups anymore, and multiplication with real numbers is restricted to non-negative ones. However, these features are exactly the ones necessary and sufficient to build a theory for convex and concave functions, i.e., for Convex Analysis. Moreover, the same type of structure forms the basis for Set Optimization: see \cite{Hamel05Habil, Hamel09SVVAN, HamelEtAl15Incoll, Loehne11Book}.

%%%
%%%New section
\section{Ordering cuts}

The order relation $\leq$ is extended from $\mathbb Q$ to $\mathcal L$ and $\mathcal U$ by defining
\begin{align*}
A \leq A' & \quad \Leftrightarrow \quad A \subseteq A', \\
B \leq B' & \quad \Leftrightarrow \quad B \supseteq B'
\end{align*}
for $A, A' \in \mathcal L$, $B, B'\in \mathcal U$. It is clearly an extension of $\leq$ on $\mathbb Q$ since for $p, p' \in \mathbb Q$ and, for example, $A = \{a \in \mathbb Q \mid a < p\}$, $A' = \{a \in \mathbb Q \mid a < p'\}$ one has $p \leq p'$ if, and only if, $A \subseteq A'$.

One easily checks that $\leq$ as defined above is a linear order on $\mathcal L$ as well as on $\mathcal U$.

\begin{remark}
\label{RemSetRelations} 
On $\mathcal L$ and $\mathcal U$ (already on $\P(\mathbb Q, -\mathbb Q_+)$ and $\P(\mathbb Q, \mathbb Q_+)$, respectively), the order relations as defined above coincide with the complete lattice extensions of set relations used in Set Optimization \cite{Hamel05Habil, Hamel09SVVAN, HamelEtAl15Incoll}. Thus, one can consider Dedekind cuts as an appropriate foundation of set relation based optimization. This might put into perspective Cantor's observation that the construction of real numbers via fundamental sequences has 'the advantage that it is most suited to the analytic calculus' (see \cite[p. 33, translation from German by AH]{EbbinghausEtAl88Book2nd}).
\end{remark}

\begin{theorem}
\label{ThmMeetSL}
The pair $(\mathcal L, \subseteq)$ is a complete lattice with least element $\vn$, greatest element $\mathbb Q$ and 
\begin{equation}
\label{EqExtendedDo}
\sup_{A \in \mathcal X} A = \bigcup_{A \in \mathcal X} A, \quad \inf_{A \in \mathcal X} A = I\of{\bigcap_{A \in \mathcal X} A}
\end{equation}
for every set $\mathcal X \subseteq \mathcal L$ where it is understood that the supremum is $\vn$ and the infimum $\mathbb Q$ if $\mathcal X = \vn$.
\end{theorem}

{\sc Proof.} For $\mathcal X \subseteq \mathcal L$, define the set
\[
\bar S := \bigcup_{A \in \mathcal X} A.
\]
Clearly, $\bar S$ satisfies (A2). It also satisfies (A3). Indeed, the set $\bar S$ is either empty, all of $\mathbb Q$ or a non-trivial subset of $\mathbb Q$ which does not have a greatest element: if $\bar a \in \bar S$ would be the greatest element, then it would also be the greatest one in $A$ for some $A \in \mathcal X$ which contradicts the fact that $A \in \mathcal L$. Thus $\bar S \in \mathcal L$. It is now a routine exercise to show that $\bar S$ is the supremum of $\mathcal X$ with respect to $\subseteq$.

Similarly, define $\bar T = I(\bigcap_{A \in \mathcal X} A)$. Since  $\bigcap_{A \in \mathcal X} A$ satisfies (A2), $\bar T$ also satisfies (A2) as well as (A3). One also has $\bar T \subseteq A$ for all $A \in \mathcal X$ and if $A' \subseteq A$ for all $A \in \mathcal X$ for some $A' \in \mathcal L$, then $A' \subseteq \bigcap_{A \in \mathcal X} A$, hence als $A' \subseteq \bar T$ since $A'$ does not include a greatest element. Hence $\bar T$ is indeed the infimum of $\mathcal X$ with respect to $\subseteq$.
\pend

\begin{remark}
\label{RemUpperSetsOrder} 
Likewise, $(\mathcal U, \supseteq)$ is a complete lattice with least element $\mathbb Q$, greatest element $\vn$  and 
\begin{equation}
\label{EqExtendedUp}
\sup_{B \in \mathcal X} B = \bigcap_{B \in \mathcal X} B, \quad \inf_{B \in \mathcal X} B = J\of{\bigcup_{B \in \mathcal X} B}
\end{equation}
for set $\mathcal X \subseteq \mathcal U$ where the supremum is $\mathbb Q$ and the infimum $\vn$ if $\mathcal X = \vn$.
\end{remark}

Constructing the real numbers via ordinary Dedekind cuts requires to show that $(\mathcal L\bs\{\mathbb Q, \vn\}, \subseteq)$ has the least upper bound property, i.e., is a conditionally complete lattice (compare, for example, \cite[p. 18, Step 3 of the proof of Theorem 1.19]{Rudin76Book}). Theorem \ref{ThmMeetSL} can be considered as a ''complete" version of this result.

%%%
%%%New section
\section{Addition for extended cuts}

The principal goal of this section is to show that is makes a difference if one introduces the addition for non-ordinary Dedekind cuts on $\mathcal L$ or on $\mathcal U$. Involving non-ordinary cuts produces two different additive structures which is an explanation for the appearance of `inf-addition' and `sup-addition' in \cite{RockafellarWets10Book3rd} and before in Moreau's work \cite{Moreau63SMUM, Moreau63RCAS, Moreau66SJL}.

One can introduce the addition on $\mathcal L$ via \eqref{EqMinkowskiAddition}, \eqref{EqEmptyAddition}: for $A, A' \in \mathcal L\bs\{\vn\}$ set
\[
A + A' = \cb{a + a' \mid a \in A, \, a' \in A'},
\]
use \eqref{EqEmptyAddition} if $\vn$ is involved and generate the cut $(A + A', \mathbb Q\bs (A+A'))$ as the sum $(A, \mathbb Q\bs A) + (A', \mathbb Q\bs A')$. This gives $(\vn, \mathbb Q) + (\mathbb Q, \vn) = (\vn, \mathbb Q)$.

Alternatively, one can introduce the addition on $\mathcal U$: for $B, B' \in \mathcal U\bs\{\vn\}$ set
\[
B + B' = \cb{b + b' \mid b \in B, \, b' \in B'},
\]
use \eqref{EqEmptyAddition} if $\vn$ is involved and generate the cut $(\mathbb Q\bs (B+B'), B + B')$ as the sum $(\mathbb Q\bs B, B) + (\mathbb Q\bs B', B')$. This gives $(\vn, \mathbb Q) + (\mathbb Q, \vn) = (\mathbb Q, \vn)$.

\begin{theorem} 
\label{ThmAdditionL}
The pair $(\mathcal L, +)$ is a commutative monoid with zero element $N_L = \{q \in \mathbb Q \mid q < 0\}$. The empty set $\vn$ is the absorbing element in $(\mathcal L, +)$. The order  $\leq$ on $\mathcal L$ is compatible with $+$ with $\vn \leq A$ for all $A \in \mathcal L$.
\end{theorem}

{\sc Proof.} This is standard for elements of $\mathcal L\bs \{\vn, \mathbb Q\}$, see \cite{Rudin76Book}. Equations \eqref{EqMinkowskiAddition}, \eqref{EqEmptyAddition} provide the necessary properties for expressions involving $\mathbb Q$, $\vn$. Especially, $\vn + N_L = N_L + \vn = \vn$ and $\mathbb Q + N_L = N_L + \mathbb Q = \mathbb Q$. Also by \eqref{EqEmptyAddition}, $\vn$ is absorbing. Finally, $A, A', C \in \mathcal L$, $A \subseteq A'$ imply $A + C \subseteq A'  + C$ and $\vn \leq A$ since this is just $\vn \subseteq A$. \pend

\begin{remark}
\label{RemAdditionU}
Likewise, the pair $(\mathcal U, +)$ is a commutative monoid with zero element $N_U = \{q \in \mathbb Q \mid 0 < q\}$. The empty set $\vn$ is the absorbing element in $(\mathcal U, +)$ and the order  $\leq$ on $\mathcal U$ is compatible with $+$ with $B \leq \vn$ for all $B \in \mathcal U$.
\end{remark}

It is important to realize that the addition of the least and the greatest element in $\mathcal L$ and $\mathcal U$ leads to different results: one has $\vn + \mathbb Q = \vn$ in both cases, but $\vn$ is the least element in $\mathcal L$ while it is the greatest in $\mathcal U$. However, both versions are consequences of the same rule, namely \eqref{EqEmptyAddition}, for the extension of the Minkowski addition for sets to $\vn$.

The addition on $\mathcal L$ generates one on $\mathcal D_L$ and the one on $\mathcal U$ a corresponding one on $\mathcal D_U$ as described before Theorem \ref{ThmAdditionL}.

\begin{example} 
\label{RemMotivateResiduation}
Defining $N_L = \{q \in \mathbb Q \mid q < 0\}$, $P_L = \{q \in \mathbb Q \mid 0 \leq q\}$ as well as  $O_L = \{q \in \mathbb Q \mid q < 1\}$, $Q_L = \{q \in \mathbb Q \mid 1 \leq q\}$ one gets the two ordinary Dedekind cuts $(N_L, P_L)$ and $(O_L, Q_L)$ which satisfy (DL3); they are generated by 0 and 1, respectively. Assume for a moment that the two non-ordinary Dedekind cuts $(\varnothing, \mathbb Q)$, $(\mathbb Q, \varnothing)$ would be additive inverses of each other, i.e., $(\varnothing, \mathbb Q) + (\mathbb Q, \varnothing) = (N_L, P_L)$ would be true. Then one would obtain the two equations
\begin{align*}
(\varnothing, \mathbb Q) + (N_L, P_L) = (\varnothing, \mathbb Q) \\
(\mathbb Q, \varnothing) + (O_L, Q_L) = (\mathbb Q, \varnothing) 
\end{align*}
by using \eqref{EqEmptyAddition} both for upper and lower sets of the cuts. However, adding the two equation would produce $(N_L, P_L) = (O_L, Q_L)$, a contradiction. Moreover, the definition $(\varnothing, \mathbb Q) + (\mathbb Q, \varnothing) = (N_L, P_L)$ is not covered by the rules for the Minkowski addition, i.e., \eqref{EqEmptyAddition}.
\end{example}

\begin{remark}
\label{RemAdditionExtReals}
The set of ordinary Dedekind cuts determines the set $\R$ of real numbers. Identifying the non-ordinary cut $(\varnothing, \mathbb Q)$ with $-\infty$ and $(\mathbb Q, \varnothing)$ with $+\infty$ and vice versa, two different versions of the addition for the extended real numbers  $\R \cup\{\pm\infty\}$ are obtained: one yields $(-\infty) \ssum (+\infty) = -\infty$, the other $(-\infty) \isum (+\infty) = +\infty$. Thus, $(\R \cup\{\pm\infty\}, \ssum)$ and  $(\R \cup\{\pm\infty\}, \isum)$ are two different commutative monoids, the first one with absorbing element $-\infty$ (and $+\infty$ quasi-absorbing in the sense of \cite{HagerHamelHeyde23ArX}), while in the second one $+\infty$ is absorbing (and $-\infty$ quasi-absorbing). Following \cite{RockafellarWets10Book3rd}, the first version is called sup-addition, denoted $\ssum$, the second one inf-addition, denoted $\isum$. Both additions of course coincide for ordinary cuts and real numbers.
\end{remark} 

\begin{theorem}
\label{ThmSupInfAdditivity}
$(\mathcal L, +, \subseteq)$ is sup-additive and $(\mathcal U, +, \supseteq)$ is inf-additive.
\end{theorem}

{\sc Proof.} One has to show for $\mathcal A, \mathcal A' \subseteq \mathcal L$:
\[
\sup\{\mathcal A + \mathcal A'\} = \sup\mathcal A + \sup\mathcal A'
\]
where all additions are in the Minkowski sense with extension \eqref{EqEmptyAddition}. Note that the addition on the left hand side is in $\mathcal P(\mathcal L)$ whereas the one on the right hand side is in $\mathcal L$.

If either $\mathcal A = \vn$ or $\mathcal A' = \vn$ or both, then one gets $\vn$ on both sides of this equation. If both $\mathcal A$ and $\mathcal A'$ are non-empty, then the result follows from 
\eqref{EqExtendedDo} since the supremum is a union. The case of $(\mathcal U, \leq)$ is parallel since the infimum is a union, see \eqref{EqExtendedUp}. \pend

A simple example shows that $(\mathcal L, +, \subseteq)$ is not inf-additive. Indeed, let $\mathcal A = \cb{A_n}_{n = 1,2, \ldots}$ with $A_n = \{q \in \mathbb Q \mid q < -1/n\}$ and $\mathcal A' = \{\mathbb Q\}$. Then $\inf \mathcal A = \vn$, $\inf \mathcal A' = \mathbb Q$ due to \eqref{EqExtendedDo} (the infimum is an intersection) which gives $\inf \mathcal A + \inf \mathcal A' = \vn$ according to \eqref{EqEmptyAddition}. On the other hand, $A_n + \mathbb Q = \mathbb Q$ for all $n$, hence $\inf (\mathcal A + \mathcal A') = \mathbb Q$, thus inf-additivity is not satisfied. 

The additions on $\mathcal L$ and $\mathcal U$ can be characterized using only proper elements, i.e., not $\vn$ and not $\mathbb Q$, via the order $\leq$. This is based on the fact that $(\mathcal L, \subseteq)$ and  $(\mathcal U, \supseteq)$ are complete lattices. The property in the following theorem was used as a definition for $\ssum$ on the extended real numbers in \cite[Def. 3.1]{HamelSchrage12JCA}. Theorem \ref{ThmSupInfAdditivity} together with \cite[Theorem 2.2]{HamelSchrage12JCA} ensures that the following constructions make sense.

\begin{theorem}
For all $A, C \in \mathcal L$ one has
\begin{align*}
A + C  & = \sup\cb{U + X \mid U, X \in \mathcal L\bs\{\vn, \mathbb Q\}, U \leq A, X \leq C} \\
	& = \bigcup \cb{U + X \mid U, X \in \mathcal L\bs\{\vn, \mathbb Q\}, U \subseteq A, X \subseteq C}.
\end{align*}
\end{theorem}

{\sc Proof.} If $A = \vn$, then there is no non-empty $U \subseteq A$, hence the supremum is taken over the empty set and yields $\vn$. The same applies if $C = \vn$. If $A = \mathbb Q$ and $C \neq \vn$, one can set $U_p = \{a \in \mathbb Q \mid a < p\}$ for $p \in \mathbb Q$ and pick $c \in C$ to get $\mathbb Q = \bigcup_{p \in \mathbb Q} \of{U_p + X}$ with $X = \{a \in \mathbb Q \mid a < c\}$, thus the supremum is $\mathbb Q$ in this case. The same applies if $A \neq \vn$ and $C = \mathbb Q$. This leaves the case $A, C \not\in \{\vn, \mathbb Q\}$: one can choose $U= A$, $X = C$ and gets the desired equality. \pend

\begin{remark}
\label{RemInfSupAddition} 
Likewise, for $B, D \in \mathcal U$ one has
\begin{align*}
B + D  & = \inf\cb{V + Y \mid V, Y \in \mathcal U\bs\{\vn, \mathbb Q\}, B \leq V, D \leq Y} \\
	& = \bigcup \cb{V + Y \mid V, Y \in \mathcal U\bs\{\vn, \mathbb Q\}, B \supseteq V, D \supseteq Y}.
\end{align*}
\end{remark}

A substitute for additive inverses is provided by the following construction.

\begin{definition} 
\label{DefPseudodifferences}
The expression
\begin{equation}
\label{EqSupRes}
A \sdif C = \sup\cb{U \in \mathcal L \mid C + U \leq A} = \bigcup\cb{U \in \mathcal L \mid C + U \subseteq A}
\end{equation}
is called the (lower) pseudodifference of $A, C \in \mathcal L$.

The expression
\begin{equation}
\label{EqSupRes}
B \idif D = \inf\cb{X \in \mathcal U \mid B \leq D + X} = \bigcup\cb{X \in \mathcal U \mid B \subseteq D + X}
\end{equation}
is called the (upper) pseudodifference of $B, D \in \mathcal U$.
\end{definition}

It should be clear from the context which version of the pseudodifference is meant. The adjectives ``lower"and ``upper" are only used if there is a danger of confusion. The pseudodifferences are actually particular cases of residuation mappings and have been used in a more general context in \cite{HamelSchrage12JCA, HamelSchrage14PJO}. 

\begin{corollary}
\label{CorSupRes} Let $A \in \mathcal L$. Then

(a) $A \sdif A = N_L$ whenever $A \neq \vn$, $A \neq \mathbb Q$,

(b) $A \sdif \vn = \mathbb Q$ and 
\[
\vn \sdif A = 
	\left\{
	\begin{array}{ccc}
	\mathbb Q & : & A = \vn \\
	\vn & : & A \neq \vn
	\end{array}
	\right.
\]

(c) $\mathbb Q \sdif A = \mathbb Q$ and 
\[
A \sdif \mathbb Q = 
	\left\{
	\begin{array}{ccc}
	\vn & : & A \neq \mathbb Q \\
	\mathbb Q & : & A = \mathbb Q 
	\end{array}
	\right.
\]
\end{corollary}

{\sc Proof.} Everything is immediate from the definition of $\sdif$. \pend

In particular, the following formulas hold true:
\begin{equation}
\label{EqInfinityDiffL}
\vn \sdif \vn = \mathbb Q \sdif \vn = \mathbb Q \sdif \mathbb Q = \mathbb Q \; \text{and} \; \vn \sdif \mathbb Q = \vn.
\end{equation}
The corresponding result for $\mathcal U$ looks as follows.

\begin{corollary}
\label{CorInfRes} Let $B \in \mathcal U$. Then

(a) $B \idif B = N_U$ whenever $B \neq \vn$, $B \neq \mathbb Q$,

(b) $B \idif \vn = \mathbb Q$ and 
\[
\vn \idif B = \left\{
	\begin{array}{ccc}
	\mathbb Q & : & B = \vn \\
	\vn & : & B \neq \vn 
	\end{array}
	\right.
\]

(c) $\mathbb Q \idif B = \mathbb Q$ and 
\[
B \sdif \mathbb Q = 
	\left\{
	\begin{array}{ccc}
	\mathbb Q & : & B = \mathbb Q \\
	\vn & : & B \neq \mathbb Q \\
	\end{array}
	\right.
\]
\end{corollary}

Strangely enough, this produces the same formulas as in the $\mathcal L$-case:
\begin{equation}
\label{EqInfinityDiffU}
\vn \idif \vn =  \mathbb Q \idif \vn = \mathbb Q \idif \mathbb Q = \mathbb Q \; \text{and} \; \vn \idif \mathbb Q = \vn.
\end{equation}

Turning this into formulas for extended real numbers, one has to identify $\vn$ with $-\infty$ and $\mathbb Q$ with $+\infty$ in $\mathcal L$-case which gives, for example, $(-\infty) \sdif (-\infty) = +\infty$ (first equation in \eqref{EqInfinityDiffL}). The identification is vice versa in the $\mathcal U$-case: $\mathbb Q$ turns into $-\infty$ and $\vn$ into $+\infty$ which yields $(+\infty) \idif (+\infty) = \infty$ (first equation in \eqref{EqInfinityDiffU}).

This procedure leads to odd looking, but well-defined identities for pseudodifferences of top and bottom elements in $\R \cup\{\pm\infty\}$ which are displayed in the table below. They can already be found in \cite{HamelSchrage12JCA} obtained without any reference to Dedekind cuts.

\bigskip
\noindent
\hspace{4.2cm} {\bf $\mathcal L$-case} \hspace{4cm} {\bf $\mathcal U$-case}
\begin{align*}
\mbox{Sum:} \quad & (+\infty) \ssum (-\infty) = -\infty & (+\infty) \isum (-\infty) = +\infty
\\[.3cm]
\mbox{Differences:} \quad
    & (+\infty) \idif (+\infty) = +\infty & \of{+\infty} \sdif \of{+\infty} = -\infty \\
    & (+\infty) \idif \of{-\infty} = +\infty & \of{+\infty} \sdif \of{-\infty} = +\infty \\
    & \of{-\infty} \idif (+\infty) = -\infty & \of{-\infty} \sdif \of{+\infty} = -\infty  \\
    & \of{-\infty} \idif \of{-\infty} = +\infty & \of{-\infty} \sdif \of{-\infty} = -\infty
\end{align*}

By the way of conclusion, not only expressions like $(+\infty) + (-\infty)$ can be given a precise algebraic meaning, but also $(+\infty) - (+\infty)$, $(-\infty) - (-\infty)$ etc. Again, there is a version based on inf-addition and another one based on sup-addition. The underlying structural property is that both $\mathcal L$ and $\mathcal U$ are residuated (complete) lattices with reversed order relation.

%%%
%%%New section
\section{Negatives and opposites} 

The element-wise multiplication by $-1$ of a set $A \in \mathcal L$ or $B \in \mathcal U$ does not produce a set in $\mathcal L$ or $\mathcal U$, respectively. One the contrary, defining
\[
(-1)X = \cb{-x \mid x \in X} 
\]
for a subset $X \subseteq \mathbb Q$ with $(-1)\vn = \vn$ one gets
\[
(-1)A \in \mathcal U \quad \text{and} \quad (-1)B \in \mathcal L
\]
for $A \in \mathcal L$, $B \in \mathcal U$. This yields the following result.

\begin{proposition}
\label{PropNegatives}
The multiplication by $-1$ is an algebraic and a reverse order isomorphism between $\mathcal L$ and $\mathcal U$, i.e.,
\[
(-1)(A + A') = (-1)A + (-1)A ' \quad \text{and} \quad (-1)(B + B') = (-1)B + (-1)B'
\]
and
\[
A \leq A'\quad \Leftrightarrow \quad (-1)A' \leq (-1)A
\]
for $A, A' \in \mathcal L$, $B, B'\in \mathcal U$.
\end{proposition}

{\sc Proof.} Straightforward. \pend

Clearly, a similar isomorphism is obtained if one replaces $-1$ by any rational number $q < 0$. In the following, $-X$ is written for $(-1)X$. Since $(-1)\vn = \vn$, $(-1)\mathbb Q = \mathbb Q$, the multiplication by $-1$ maps the least and the greatest element of $\mathcal L$ onto the greatest and the least in $\mathcal U$, respectively, and vice versa.

As usual, one can express inverses with respect to addition using the multiplication by $-1$. If $A \not\in \{\vn, \mathbb Q\}$, then 
\[
A + I\of{-\mathbb Q\bs A} = N_L \quad \text{and} \quad B + J\of{-\mathbb Q\bs B} = N_U.
\]
Indeed, $a \in A$ and $x \in I\of{-\mathbb Q\bs A}$ means that $x \in -\mathbb Q\bs A$ but is not the least element in that set. Moreover, this implies $-x \in \mathbb Q\bs A$, hence $a < -x$ and $a + x < 0$, so $x + a \in N_L$. Finally, if $q < 0$, then take $a \in A$ such that $-a \in I\of{-\mathbb Q\bs A}$ (always exists since $A \in \mathcal L\bs\{\vn, \mathbb Q\}$) and get $(a + q) + (-a) = q < 0$, hence every $q \in N_L$ can be represented as a sum of an element of $A$ and one of $I\of{-\mathbb Q\bs A}$. Similar for the $B$-case.

In the following, the symbols
\[
A^* = I\of{-\mathbb Q\bs A} \quad\text{and}\quad B^* = J\of{-\mathbb Q\bs B} 
\]
are used to denote the additive inverses of $A \in \mathcal L\bs\{\vn, \mathbb Q\}$, $B \in \mathcal U\bs\{\mathbb Q, \vn\}$. However, these expressions are also defined for $A , B \in \{\vn, \mathbb Q\}$ but are not additive inverses anymore: $I\of{-\mathbb Q\bs \mathbb Q} = \vn = J\of{-\mathbb Q\bs \mathbb Q} $, $ I\of{-\mathbb Q\bs \vn}  = \mathbb Q = J\of{-\mathbb Q\bs \vn}$.

%%%
%%%New section
\section{Multiplication} 

Denote 
\[
\mathcal L_+ = \cb{S \in \mathcal L \mid  N_L \subseteq S}  \quad \text{and} \quad \mathcal U_+ = \cb{T \in \mathcal U \mid  N_U \supseteq T}
\]
and set $\mathcal L_+' = \mathcal L_+ \bs\{\mathbb Q\}$, $ \mathcal U_+' =  \mathcal U_+ \bs \{\vn\}$. Note that $\vn \not\in \mathcal L_+$, $\mathbb Q \not\in \mathcal U_+$.
 
The products of two elements $S \in \mathcal L_+'$, $A \in \mathcal L_+$ and $T \in \mathcal U_+'$, $B \in \mathcal U_+$ are defined by
\begin{align*}
S \cdot A & = \cb{sa \mid s \in S, a \in A \; \text{and} \; (s \geq 0, a \geq 0 \; \text{or} \; sa < 0)} \\
T \cdot B & = \cb{tb \mid t \in T, b \in B}
\end{align*}
with $N_U \cdot \vn = N_U$ and $T \cdot \vn = \vn$ for $T \in \mathcal U_+'\bs\{N_U\}$. The presence of negative numbers in $S$, $A$ makes the first case more cumbersome. With this definition, one gets
\[
N_L \cdot A  = N_L, \; O_L \cdot A = A, \; N_U \cdot B = N_U, \; O_U \cdot B = B
\]
for all $A \in \mathcal L_+$, $B \in \mathcal U_+$.

The multiplication is extended to $A \in \mathcal L$ by setting
\[
S \cdot A = I\of{-\mathbb Q\bs (S \cdot I(-\mathbb Q\bs A))} 
\]
whenever $A \not\in \mathcal L_+$ with $N_L \cdot \vn = N_L$ and $S \cdot \vn = \vn$ for $S \in \mathcal L_+'\bs\{N_L\}$. Likewise,
\[
T \cdot B = J\of{-\mathbb Q\bs (T \cdot J(-\mathbb Q\bs B))}
\]
for $B \not\in \mathcal U_+$.

If $A \in \mathcal L\bs\{\vn, \mathbb Q\}$, $S \in \mathcal L_+'$ and $B \in \mathcal U\bs\{\mathbb Q, \vn\}$, $T \in \mathcal U_+'$ then
\[
S \cdot A = (S \cdot A^*)^* \quad\text{and}\quad T \cdot B = (T \cdot B^*)^*
\]
which corresponds to known definitions such as \cite[p. 20, Step 7]{Rudin76Book}.

\begin{theorem}
\label{ThmMultiplicationLower}
If $S, S_1, S_2 \in \mathcal L_+'$ and $A, A_1, A_2 \in \mathcal L$, then one has

(a) $S \cdot \of{A_1 \ssum A_2} = (S \cdot A_1) \ssum (S \cdot A_2)$,

(b) $S_1 \cdot \of{S_2 \cdot A} = \of{S_1\cdot S_2} \cdot A$,

(c) $O_L \cdot A = A$,

(d) $N_L \cdot A = N_L$.

Finally, if $A_1 \leq A_2$, then $S \cdot A_1 \leq S \cdot A_2$.
\end{theorem}

{\sc Proof.} This result is well-established if all involved cuts are ordinary. It remains to be checked if the arithmetic rules also apply to non-ordinary cuts which is a straightforward exercise. $\square$

\begin{corollary}
Both $(\R\cup\{\pm\infty\}, \ssum, \cdot)$ and $(\R\cup\{\pm\infty\}, \isum, \cdot)$ are conlinear spaces as defined in \cite{Hamel05Habil} where $\cdot$ denotes the multiplication with non-negative real numbers.
\end{corollary}

%%%
%%%New section
\section{Applications}

%%%
%%%New subsection
\subsection{Convex functions}

Rockafellar \cite[p. 23]{Rockafellar70Book} defined convexity of extended real-valued functions via their epigraphs, i.e., a function $f \colon \R^n \to \R\cup\{\pm\infty\}$ is convex if the set
\[
\epi f = \cb{(x, r) \in \R^n \times \R \mid f(x) \leq r}
\] 
is convex. Jensen's inequality as equivalent condition can only be used if expressions like $(+\infty) + (-\infty)$ are avoided which is discussed in \cite[p. 24]{Rockafellar70Book}. Passing to inf-addition, the difficulty disappears.

\begin{theorem}
\label{ThmConvexity} Let $X$ be a non-trivial real linear space. The following two conditions are equivalent for a function $f \colon X \to \R\cup\{\pm\infty\}$:

(a) $\epi f = \cb{(x, r) \in X \times \R \mid f(x) \leq r}$ is a convex subset of $X \times \R$,

(b) one has
\begin{equation}
\label{EqJensenIneq}
f(\alpha x + (1- \alpha) y) \leq \alpha f(x) \isum (1-\alpha) f(y) 
\end{equation}
whenever $x, y \in X$ and $\alpha \in (0,1)$.
\end{theorem}

{\sc Proof.} (a)$\Rightarrow$(b): The result is well-established for $f(x), f(y) \in \R\cup\{+\infty\}$.

If $f(x) = f(y) = -\infty$, then $(x, r), (y, r) \in \epi f$ for all $r \in \R$, hence by assumption
\[
\forall r \in \R \colon (\alpha x + (1-\alpha) y, r) \in \epi f.
\]
Therefore $f(\alpha x + (1-\alpha) y) = -\infty$ which means that \eqref{EqJensenIneq} holds true. Finally assume $f(x) = -\infty$. If $f(y) =  +\infty$, then Jensen's inequality holds true since $(-\infty) \isum (+\infty) = +\infty$. If $f(y) \in \R$, then $(y, f(y)) \in \epi f$ and $(x, r) \in \epi f$ for all $r \in \R$, hence by assumption
\[
\forall r \in \R \colon (\alpha x + (1-\alpha) y, \alpha r + (1-\alpha)f(y)) \in \epi f
\]
implying $f(\alpha x + (1-\alpha) y) = -\infty$ and Jensen's inequality is satisfied also in this case.

(b)$\Rightarrow$(a): Assume $(x, r), (y, s) \in \epi f$. Using \eqref{EqJensenIneq} one gets
\[
f(\alpha x + (1- \alpha) y) \leq \alpha f(x) \isum (1-\alpha) f(y) \leq \alpha s + (1-\alpha) r,
\]
hence $(\alpha x + (1- \alpha) y, \alpha s + (1-\alpha) r) \in \epi f$ which is convexity of $\epi f$. \pend

The equivalence in Theorem \ref{ThmConvexity} is no longer true if the inf-addition is replaced by the sup-addition $\ssum$. An example is the function $f \colon \R \to \R\cup\{\pm\infty\}$ defined by (see \cite[p. 41]{RockafellarWets10Book3rd})
\[
f(x) = \left\{
\begin{array}{ccc}
+\infty & : & x < 0 \\
0 & : & x = 0 \\
-\infty & : & x > 0 
\end{array}
\right..
\]
It is convex, \eqref{EqJensenIneq} is satisfied, but one has $0 = f(0) > 1/2f(-1) \ssum 1/2f(1) = (-\infty) \ssum (+\infty) = -\infty$.

\begin{remark}
The situation is reversed if concavity is the focus. The equivalence in Theorem \ref{ThmConvexity} turns to the equivalence between the convexity of the hypograph and Jensen's inequality with reversed inequality sign and $\ssum$ instead of $\isum$. Thus, the theory is perfectly symmetric despite 'there is no single, symmetric way of handling $\infty - \infty$' \cite[p. 15]{RockafellarWets10Book3rd}.
\end{remark}

Two functional operations, widely used for convex functions, crucially depend on addition: the pointwise addition of functions and the infimal convolutions. In \cite[p. 33, 34]{Rockafellar70Book}, it is carefully explained why the corresponding theorems are only formulated for proper convex functions, namely 'for the sake of avoiding $\infty-\infty$.'

\begin{theorem}
\label{ThmFunctOperations}
Let $f_1, f_2 \colon X \to \R\cup\{\pm\infty\}$ be two convex functions. Then

(a) $f(x) = (f_1 \isum f_2)(x) := f_1(x) \isum f_2(x)$ is convex,

(b) $f(x) = \inf\{f_1(x_1) \isum f_2(x_2) \mid x_1 + x_2 = x\}$ is convex.
\end{theorem}

{\sc Proof.} (a) Straightforward using \eqref{EqJensenIneq}. (b) One can use Theorem \ref{ThmConvexity} as follows. Take $(x, r), (y, s) \in \epi f$. Fix $\eps > 0$. Then, there are $x^\eps_1, x^\eps_2, y^\eps_1, y^\eps_2 \in X$ satisfying $x^\eps_1+ x^\eps_2 = x$, $y^\eps_1+ y^\eps_2 = y$ and
\begin{align*}
f_1(x^\eps_1) \isum f_2(x^\eps_2) & \leq r + \eps, \\
f_1(y^\eps_1) \isum f_2(y^\eps_2) & \leq s + \eps.
\end{align*}
In particular, $f_1(x^\eps_1), f_2(x^\eps_2), f_1(y^\eps_1), f_2(y^\eps_2) \neq +\infty$. The two inequalities above produce
\[
\alpha f_1(x^\eps_1) \isum (1-\alpha) f_1(y^\eps_1) \isum \alpha f_2(y^\eps_1) \isum (1-\alpha) f_1(y^\eps_2) \leq \alpha r + (1-\alpha)s + \eps.
\]
Jensen's inequality \eqref{EqJensenIneq} further yields
\[
 f_1(\alpha x^\eps_1 + (1-\alpha) y^\eps_1) \isum  f_2(\alpha y^\eps_1 + (1-\alpha) y^\eps_2) \leq \alpha r + (1-\alpha)s + \eps.
\]
Since this is true for all $\eps > 0$, one gets the same inequality with $\eps = 0$ and consequently
$\alpha (x,r) + (1-\alpha)(y,s) \in \epi f$. 
\pend

The reader may verify that both statements are not true if the inf-addition $\isum$ is replaced by the sup-addition $\ssum$. 

\begin{remark}
One may also compare (b) of Theorem \ref{ThmFunctOperations} with (ix) of  Theorem 2.1.3 in the benchmark text \cite{Zalinescu02Book}: therein, the two functions $f_1, f_2$ are assumed to be proper.
\end{remark}

Thus, the consideration of extended real-valued functions together with proper arithmetic rules admits to deal with the proper and improper case in a unified way.

%%%
%%%New subsection
\subsection{Set-valued convex function}

In this subsection, one issue related to set-valued convex functions and set optimization is indicated. The reader is referred to the literature such as \cite{Loehne11Book}, \cite{Hamel09SVVAN, HamelEtAl15Incoll, HamelSchrage12JCA, HamelSchrage14PJO},  for more details and applications.

Let $C \subseteq \R^d$ be a closed convex cone with $C \not\in \{\vn, \R^d\}$. A function $f \colon X \to \mathcal P(\R^d)$ is called $C$-convex iff
\begin{equation}
\label{EqJensenIneqSet}
f(\alpha x + (1- \alpha) y) + C \supseteq \alpha f(x) + (1-\alpha) f(y) 
\end{equation}
whenever $x, y \in X$ and $\alpha \in (0,1)$. Extending $f$ to $f^\triup$ defined by $f^\triup(x) = \cl(f(x) + C)$ if $f(x) \neq \vn$ and $f^\triup(x) = \vn$ otherwise one gets a function $f^\triup$ which maps $X$ into 
\[
\mathcal G(\R^d, C) = \cb{A \subseteq \R^d \mid A = \cl\co(A + C)}.
\]
Using standard separation arguments one can show that such functions can be characterized by families of scalar functions, namely
\[
\vp^\triup_{f, w}(x) = \inf_{z \in f^\triup(x)} w^\top x
\]
where $w$ runs through the dual cone $C^+ = \{w \in \R^d \mid \forall z \in C \colon w^\top z \geq 0\}$, see \cite[formula (34)]{HamelSchrage12JCA} (with a different notation) or \cite[formula (4.17)]{HamelEtAl15Incoll}.

The following example \cite[Example 5.25]{HamelSchrage12JCA} shows that scalarizations of type $\vp^\triup_{f, w}$ can be improper even for well-behaved set-valued functions.

\begin{example}
\label{ExImproperScalarization}
Let $X = \R$, $d = 2$, $v =  (0, 1)^\top$ and $C = \{sv \mid s \geq 0\}$ and consider the function
\[
f(x) = f^\triup(x) =
\left\{
	\begin{array}{ccc}
	H^+(v) & : & x > 0 \\
	C & : & x = 0 \\
	\vn & \colon & x < 0
	\end{array}
	\right.
\]
where $H^+(v) = \{z \in \R^2 \mid v^\top z = z_2 \geq 0\}$. This function is $C$-convex and $C$-proper in the sense that $\dom f = \{x \in \R \mid f(x) \neq \vn\} \neq \vn$ and $f(x) \neq \R^2$ for all $x \in \R$ (see \cite{Hamel09SVVAN}). However, the scalarization $\vp^\triup_{f, w}$ is proper in the traditional sense, i.e., $\dom \vp^\triup_{f, w} \neq \vn$ and $\vp^\triup_{f, w}(x) \neq -\infty$ for all $x \in \R$, only if $w$ is collinear to $v$. For example, for $w = (1, 1)^\top$ one gets
\[
\vp^\triup_{f, w}(x) = 
\left\{
	\begin{array}{ccc}
	-\infty & : & x > 0 \\
	0 & : & x = 0 \\
	+\infty & \colon & x < 0
	\end{array}
	\right.
\]
Moreover, one cannot characterize $f$ only with its proper scalarizations since $f(0) = C \subsetneq H^+(v)$ is true.
\end{example}

Thus, the consideration of the extended real numbers $\R\cup\{\pm\infty\}$ together with proper arithmetic rules admits to deal with proper and improper scalarization at the same time avoiding to impose particular restrictions to set-valued functions ensuring the properness of their scalarizations. Moreover, this also leads to appropriate conlinear spaces of scalarizing functions, basically support functions of values of a set-valued function, as developed in \cite[Section 4.2]{HamelEtAl15Incoll}.

\medskip\noindent
{\bf Acknowledgement.} No external funding was used during the preparation of this work. There are no competing interests.

\medskip\noindent
{\bf Data declaration.} There are no additional data related to this work.

%%%Bibo

\end{document}